\documentclass[12pt, leqno]{article}

\usepackage{amsmath, amsthm, amssymb, amsfonts, enumerate}
\usepackage[applemac]{inputenc} 
 \usepackage[french] {babel}

\def \R{I\!\!R}

\def \N{I\!\!N}

\newtheorem{thm}{Theorem}[section]
\newtheorem{cor}[thm]{Corollary}
\newtheorem{lem}[thm]{Lemma}
\newtheorem{pro}[thm]{Proposition}
\newtheorem{defi}[thm]{Definition}
\newtheorem{rem}[thm]{Remark}

\def \Max{ {\rm Max} }

\newcommand{\eps}{\varepsilon}

\newcommand{\be}{\begin{equation}}
\newcommand{\ee}{\end{equation}}

\def\QED{\begin{flushright} {\bf QED } \end{flushright}}

\newcommand{\dok}{\noindent {\bf Proof} : }

\def\abstract{\begin{center} \small\bf Abstract\end{center}\small}

\title{On the existence of a limit  value in some non expansive optimal
control problems.}

\author{Marc Quincampoix\thanks{Laboratoire de Math\'ematiques, UMR6205, Universit\'e de
Bretagne Occidentale, 6 Avenue Le Gorgeu, 29200 Brest, France.
Marc.Quincampoix@univ-brest.fr} , J\'er\^ome Renault\thanks {GIS ``Decision
Sciences" X-HEC-ENSAE, CMAP
and Department of Economics, Ecole Polytechnique,
91128 Palaiseau Cedex, France.  jerome.renault@polytechnique.edu}}

\date{July 4 th,  2009}

\begin{document}

\maketitle
\begin{abstract}
We investigate a limit value of an optimal control problem when
the horizon converges to infinity. For this aim, we suppose
suitable nonexpansive-like assumptions which does not imply  that
the limit is independent of the initial state as it is usually
done in the literature.
\end{abstract}

\section{Introduction}

We consider the following optimal  control denoted $\Gamma_t(y_0) $:
\begin{equation} \label{eq1} V_t(y_0) := \inf_{u \in {\cal U}} \;  \frac{1}{t} \int_{s=0}^t h(y(s,u, y_0),u(s))ds,\end{equation}
where $s\mapsto y(s,u, y_0) $ denotes the solution to
\begin{equation} \label{eq2} y'(s)=g(y(s),u(s)), \;\;\;y(0)=y_0.  \end{equation}
Here ${\cal U}$ is the set of measurable controls  from $\R_+$ to
a given non empty     metric space $U$. Throughout
the paper, we will suppose Lipschitz regularity of $g: \R ^d
\times U \to \R^d $ which  implies that for a given control $
u$ in ${\cal U}$ and a given initial condition $y_0$, equation (\ref{eq2}) has a unique absolutely
continuous solution.

The main goal of the paper consists in studying  the asymptotic
behaviour of $V_t (y_0)$  when $t$ tends to  $ \infty$. This
problem has been considered in several papers  (cf for instance in
\cite{Arisawa, Bettiol, Carda}) by approaches ensuring that the
limit of $V_t (y_0)$ is independent of $y_0$. In the present paper
we exhibit several examples where the limit exists and depends of
$y_0$. Our aim is to obtain a general result which contains in
particular the more  easy to state  following result, where
throughout the paper,  $< \cdot, \cdot>$ stands for the canonical
scalar product  and $B$ is  the associated closed unit ball..
\begin{pro} \label{proint}  Assume  that $g$ is Lipschitz, that there exists a compact set $N$ which is - forward - invariant by the control system (\ref{eq2}) and that $h$ is a continuous function which does not depend on $u$. Assume moreover that:
\be \label{nonex} \forall (y_1,y_2)  \in N^2, \; \sup_{u \in U}
\inf_{v \in U} <y_1-y_2,g(y_1,u)-g(y_2,v)>\leq 0.\ee

Then problem (\ref{eq1}) has a value when $t$ converges to $+\infty$  i.e.
there exists $V(y_0):=  \lim _{t \to + \infty } V_t(y_0)$. \end{pro}
Condition (\ref{nonex}) means a non expansive  property of the control
 system, while the condition
 $$\forall (y_1,y_2)  \in N^2, \; \sup_{u \in U} \inf_{v \in U} <y_1-y_2,g(y_1,u)-g(y_2,v)>\leq -C \|y_1 - y_2\|^2$$ expresses a dissipativity property of the control system. The above dissipativity condition does imply that the limit is independent of  $y_0$ (cf \cite{artstein}).

The value function (\ref{eq1}) can also be characterized through  - viscosity - solution of  a suitable
 Hamilton-Jacobi equation. In several articles initiated by the pioneering work \cite{Lions} the limit of $V_t(y_0)$ is obtained by ``passing to the limit" on the Hamilton-Jacobi equation.  This required coercivity properties of the Hamiltonian which could be implied by controlability and/or dissipativity of the control system but which are not valid in the nonexpansive case (\ref{nonex}). Moreover the PDE
approach is out of the scope of the - long enough -  present article.

 \begin{defi} $\;$
The problem $\Gamma(y_0):=(\Gamma_t(y_0))_{t >0}$ has a {\it limit value} if  $\lim_{t \to \infty} V_t(y_0)$ exists. Whenever it exists, we denote this limit by $V(y_0)$.
\end{defi}

Our main aim consists in giving one sufficient condition
ensuring the existence of the limit value. As a particular case of our main result we obtain proposition (\ref{proint}).

It is also of  interest to know if approximate optimal
controls for the value $V_t(y_0)$ are still approximate optimal
controls for the limit value. This leads us to the following definition.

 \begin{defi} $\;$
The problem $\Gamma(y_0)$  has a   {\it uniform value} if it has a limit value $V(y_0)$ and if:
  $$\forall \varepsilon >0, \exists u \in {\cal U},  \exists t_0,
    \forall t \geq t_0, \frac{1}{t} \int_{s=0}^t h(y (s,u,y_0),u(s))ds\leq  V(y_0)+\varepsilon.$$
\end{defi}

 Whenever the uniform value exists, the controller can act (approximately) optimally independently of the time horizon. On the contrary, if the limit value exists but the uniform value does not, he   really needs to know the time horizon before choosing a control. We will prove that our results do imply the existence of a uniform value. We will be inspired by a recent work in the discrete time case \cite{Renault}.

Let us explain now, how the paper is organized. The second section
contains some preliminaries and discussions of limit behaviors in
examples.
 In the third section, we state and prove our main result for the existence of the uniform value.

\section{Preliminaries}

We now consider the optimal control problems $(\Gamma _t(y_0))_t$
described by (\ref{eq1}) and (\ref{eq2}).

\subsection{Assumptions and Notations}

We now describe the assumptions made on $g$ and $h$.  \be
\label{HYPO} \left\{
\begin{array}{l}
\mbox{The function  $h: \R^d\times U\longrightarrow \R $ is
 measurable and bounded} \\
 \mbox{The function  $g: \R^d\times U\longrightarrow \R^d$ is
 measurable} \\
 \mbox{$\exists L
\geq 0,\forall (y,y') \in \R ^{2d} , \forall u \in U, \;
\|g(y,u)-g(y',u)\|\leq L \|y-y'\|$} \\ \mbox{$\exists a>0, \forall
(y,u) \in \R ^d \times U , \;  \|g(y,u)\|\leq a(1+\|y\|)$ }
\end{array}
\right. \ee

  With these hypotheses, given $u$ in ${\cal U}$ equation (\ref{eq2})
  has a unique absolutely continuous solution $y(\cdot, u,y_0): \R_+ \rightarrow \R^d$.

  Since $h$ is bounded, we will assume without loss of generality
from now on that $h$ takes values in $[0,1]$.

We denote by $G(y_0):=\{y(t, u,y_0), t\geq 0, u\in {\cal U}\}$ the
reachable set (i.e. the set of states that can be reached starting
from $y_0$).

We denote the average cost induced by $u$ between time 0 and time $t$ by:
  $$\gamma_{t}(y_0,u)=\frac{1}{t} \int_{0}^{t} h(y(s,
  u,y_0),u(s))ds$$
The corresponding Value function satisfies
 $V_{t}(y_0)= \inf_{u \in {\cal U}} \gamma_{t}(y_0,u).$

%

 \subsection{Examples}

 We present here basic examples. In all these examples, the cost
  $h(y,u)$ only depends on the state $y$. We will prove later that the uniform value exists in examples 2, 3 and 4. \\

 $\bullet$  Example 1: here $y$ lies in $\R^2$ seen as the complex plane,
 there is no control and the dynamic is given by $g(y,u)= i \; y$, where $i^2=-1$.  We clearly  have:
  $$V_t(y_0) \xrightarrow[t \to \infty]{} \frac{1}{2 \pi |y_0|} \int_{|z|=|y_0|} h(z)dz,$$
 and since there is no control,  the value is uniform. \\

$\bullet$ Example 2: in the complex plane again, but now $g(y,u)= i \; y \; u $, where $u\in U$  a given bounded subset of $\R$, and $h$ is continuous in $y$.  \\

 $\bullet$ Example 3: $g(y,u)=-y+u$, where $u \in U$   a given bounded subset of $\R^d$, and $h$ is continuous in $y$.\\

$\bullet$ Example 4:  in $\R^2$. The initial state is $y_0=(0,0)$
and  the control set is  $U=[0,1]$. For a state $y=(y_1,y_2)$ and
a control $u$, the dynamic is given by $y'(s) =g(y(s) ,u(s))=
\left(
\begin{array}{c}
u(s)(1-y_1(s)) \\
 u^2(s)(1-y_1(s)) \\
 \end{array}\right)$, and the cost is  $h(y)=1 -y_1(1-y_2)$.
 Notice that for any control, $y_1 '(s)  \geq  y_2 '(s) \geq  0$,
 and  thus $y_2(t)\leq y_1(t)$ for each $t\geq 0$. One can easily
 observe that $G(y_0) \subset [0,1]^2$.

 If one uses the constant control $u=\varepsilon >0$, we obtain
 $y_1(t)=1-{\rm exp}(-\varepsilon t)$ and $y_2(t)= \varepsilon y_1(t)$.
  So we have $V_t(y_0) \xrightarrow[t \to \infty]{} 0.$

 More generally, if the initial state is $y=(y_1,y_2)\in [0,1]^2$,
 by choosing a constant control $u=\varepsilon>0$ small,
 one can show that the limit value exists and $\lim_{t \to \infty} V_t(y)=y_2$.

 Notice that there is no hope here   to use  an ergodic property,
 because $$\{y \in [0,1]^2, \lim_{t \to \infty} V_t(y)= \lim_{t \to \infty} V_t(y_0)\}= [0,1]\times \{0\},$$and starting
  from $y_0$ it is possible to reach no point in  $(0,1]\times \{0\}$.\\

$\bullet$  {Example 5:} in $\R^2$, $y_0= (0,0)$, control set
$U=[0,1]$, $g(y,u)=(y_2,u)$, and $h(y_1,y_2)=0$ if $y_1\in [1,2]$,
$=1$ otherwise.

We have  $u(s)=y_2'(s) = y_1 ''(s) $, hence we may think of the
control $u$ as the acceleration, $y_2$ as the speed and $y_1$ as
the position of some mobile. If $u=\varepsilon$ constant, then
$y_2(t)=\sqrt{2 \varepsilon y_1(t)}$ $\forall t\geq 0$.

We have $u \geq 0$, hence the speed cannot decrease. Consequently,
the time interval where $y_1(t)\in [1,2]$ cannot be longer  than
the time interval where $y_1(t)\in [0,1)$, and we have
$V_T(y_0)\geq 1/2$ for each $T$.

 One can prove that $V_T(y_0) \xrightarrow[T \to \infty]{} 1/2$ by considering the following controls:
 choose $\hat{t}$ in $(0,T)$  such that $(2 / \hat{t})  + (\hat{t}      / 2 )=T$,
  make a full acceleration up to $\hat{t}$ and completely stop accelerating after:
  $u(t)=1$ for $t <{\hat{t}}$, and $u(t)=0$ for $t\geq \hat{t}$.

  Consequently the limit value exists and is $1/2$.
  However, for any control $u$ in ${\cal U}$, we either have
  $y(t, u,y_0)=y_0$ for all $t$, or $y_1(t, u,y_0) \xrightarrow[t \to \infty]{} + \infty$.
   So in any case we have  $\frac{1}{t} \int_{0}^{t} h(y(s, u,y_0),u(s))ds \xrightarrow[t \to \infty]{} 1$.
   The uniform value does not exist here, although the dynamic is very regular.

 \section{Existence results for the uniform value}

\subsection{A technical Lemma} Let us define $V^-(y_0):=\liminf_{t
\to + \infty } V_t(y_0) $ and $V^+(y_0):=\limsup_{ t \to + \infty
} V_t(y_0)$.
 Adding a parameter $m\geq 0$,  we will more generally consider the costs between  time $m$ and time $m+t$:
   $$\gamma_{m,t}(y_0,u)=\frac{1}{t} \int_{m}^{m+t} h(y(s, u,y_0),u(s))ds,$$
    \noindent and the value of the problem where the time interval $[0,m]$ can be devoted to reach a good initial state, is denoted by:
   $$V_{m,t}(y_0)= \inf_{u \in {\cal U}} \gamma_{m,t}(y_0, u) .$$
\noindent Of course $\gamma_{t}(y_0,u)= \gamma_{0,t}(y_0,u)$ and
$V_t({y_0})=V_{0,t}(y_0)$.

 \begin{lem} \label{lem1} For every  $m_0$ in $\R_+$, we have:
$$ \sup _{t>0 } \inf _{m\leq m_0}  V_{m,t}( y_0) \geq V ^+ (y_0) \geq V ^- (y_0)
\geq \sup _{t>0 } \inf _{m \geq 0 } V_{m,t} (y_0).$$ \end{lem}

\noindent{\bf Proof:}   We first prove $\sup _{t>0 } \inf _{m\leq
m_0}  V_{m,t}( y_0) \geq V ^+ (y_0)$. Suppose by contradiction
that it is false. So there exists $ \eps >0 $ such that for any
$t>0$ we have $  \inf _{m\leq m_0}  V_{m,t}( y_0) \leq V ^+ (y_0)
 - \eps $ . Hence  for any $t>0$ there exists $m \leq m_0$ with $
V_{m,t}( y_0) \leq V ^+ (y_0) - (\eps /2) $.  Now observe that
$$ \displaystyle{V_{m,t}( y_0)
= \inf _{u} \frac{1}{t}\int _m ^{m+t}h(y(s, u,y_0),u(s))ds =
\frac{1}{t} \inf_{u} \{\int _0 ^{m_0 +t }h(y(s, u,y_0),u(s))ds }$$
$$ \displaystyle{  - \int _{m+t}  ^{m_0+t}h(y(s,
u,y_0),u(s))ds - \int _0 ^{m}h(y(s, u,y_0),u(s))ds \} \geq
\frac{m_0+t}{t} V_{m_0+t} (y_0) -2 \frac{m_0}{t}.  }$$ Hence
$$\frac{m_0+t}{t} V_{m_0+t} (y_0) -2 \frac{m_0}{t} \leq V ^
+(y_0) -(\eps /2).$$ Passing to the limsup when  $t$ goes to $+
\infty $ we obtain a contradiction.

We now prove $V ^- (y_0) \geq \sup _{t>0 } \inf _{m \leq 0 }
V_{m,t} (y_0).$
  Assume
on the contrary that it is false. Then there exists  $ \eps >0 $
and $t>0$    such that $V ^-
(y_0) + \eps  \leq  \inf _{m \leq 0 } V_{m,t} (y_0).$ So  for any
$m \geq 0$,  we have $V ^- (y_0) + \eps  \leq   V_{m,t} (y_0).$
 We will obtain a contradiction
by concatenating trajectories. Take $T>0$, and write $T=lt+r$,
with $l$ in $\N$ and $r$ in $[0,t)$. For any control $u$ in ${\cal
U}$, we have: $T \gamma_T(y_0,u)$ $=$ $t \gamma_{0,t}(y_0,u) $ $+$
$  t \gamma_{t,t}(y_0,u)$ $ +$ $ ...$ $ +$ $ t
\gamma_{(l-1)t,t}(y_0,u)$ $ +$ $ r \gamma_{lt,r}(y_0,u)$ $\geq l t  (V ^- (y_0) + \eps )$. Hence
$$\displaystyle{ \gamma_T(y_0,u)
\geq \frac{T-r}{T} (V ^- (y_0) + \eps ).}$$ So for $T$ large
enough  we have $V_T(y_0)\geq V^-(y_0) +  \varepsilon/2$, hence a
contradiction by taking the liminf when $T \to \infty$ . \hfill
$\Box$

 \vspace{0,5cm}

   \noindent   {\bf  Remark:}  it
   is also easy to show that for each $t_0\geq 0$,
    we have $ \inf_{m\geq 0}  \sup_{t>t_0}  V_{m,t}(y_0)\geq V^+(y_0).$ \\

The following quantity will play a great role in the sequel.

\begin{defi}  $$V^*(y_0)= \sup _{t> 0}\; \inf _{m\geq 0} \; V_{m,t}(y_0).$$ \end{defi}

\subsection{Main results}
Let us state the first version of our main result (which clearly
implies Proposition \ref{proint} stated in the introduction)
\begin{pro} \label{thm1}  Assume  that (\ref{HYPO}) holds true and furthermore:

\noindent (H'1) $h(y,u)=h(y)$ only depends on the state,  and  is
continuous on $\R^d$.

\noindent (H'2)   $G(y_0)$ is bounded,

\noindent (H'3) $\forall (y_1,y_2)  \in G(y_0) ^2 $, $\;\; \sup_{u
\in U} \inf_{v \in U} <y_1-y_2,g(y_1,u)-g(y_2,v)>\leq 0.$

\noindent Then the problem $\Gamma(y_0)$ has a {\it limit value}
which is $V^*(y_0)$, i.e. $V_t(y_0) \xrightarrow[t \to +\infty]{}
V^*(y_0)$. The convergence of $(V_t)_t$ to $V^*$  is uniform over
$G(y_0)$, and we have

\noindent $V^*(y_0)= \sup_{t \geq 1}\; \inf_{m\geq 0} \;
V_{m,t}(y_0)=
 \inf_{m\geq 0}\sup_{t \geq 1} \; V_{m,t}(y_0)$ $=$ $\lim_{m \to \infty, t \to \infty}
 V_{m,t}(y_0)$. Moreover the value of $\Gamma(y_0)$  is uniform. \end{pro}

Condition ({\it H'3}) can be used to show that (cf Proposition
\ref{pro1}): $\forall (y_1,y_2) \in G(y_0) ^2 $ , $\forall
\varepsilon>0$, $\forall T \geq 0$, $\forall u \in {\cal U}$,
$\exists v \in {\cal U}$ s.t.: $\forall t \in [0,T], \;\; \|
y(t,u,y_1) -  y(t,v,y_2) \|\leq \|y_1 -y_2\|+ \varepsilon.$
Proposition \ref{thm1} can be applied to the previous examples 1,
2 and 3, but not to example 4.
Notice that in example 5, we have $V^*(y_0)=0 < 1/2= \lim_t V_t(y_0).$ \\

We will prove the following generalization of Proposition
\ref{thm1}. We put $Z=G(y_0)$, and denote by $\bar{Z}$ its closure
in $\R^d$.


\begin{thm} \label{thm2} Suppose that (\ref{HYPO}) holds true and furthermore assume  that

\noindent (H1)  $h$ is  uniformly continuous in $y$ on $\bar{Z} $
uniformly in $u$.
 And for each $y$ in $\bar{Z}$,  either $h$ does not depend on $u$ or the set $\{(g(y,u), h(y,u)) \in \R^d \times [0,1] , \;
 u \in U\}$ is closed.

\noindent(H2): There exist a continuous  function   $\Delta:\R^d \times
\R^d \longrightarrow \R_+$,
vanishing on the diagonal ($\Delta(y,y)=0$ for each $y$) and
symmetric ($\Delta(y_1,y_2)=\Delta(y_2,y_1)$
  for all $y_1$ and $y_2$), and a function $\hat{\alpha}:  \R_+ \longrightarrow \R_+$    s.t.  $ \hat{\alpha}(t) \xrightarrow[t \to 0]{} 0$ satisfying:

 $a)$  For every sequence $(z_n)_n$ with values in ${Z}$ and every $\varepsilon>0$,
 one can  find $n$ such that $\liminf_p \Delta(z_n,z_p)\leq \varepsilon$.

$b)$ $\forall (y_1,y_2)  \in \bar{Z}^2 ,   \; \forall u \in U$,
$\exists v \in U$ such that

 \centerline{$D \uparrow \Delta(y_1,y_2)(g(y_1,u),g(y_2,v))\leq 0$ and $h(y_2,v)-h(y_1,u)
  \leq \hat{\alpha}(\Delta(y_1,y_2)).$}

\vspace{0,3cm}

Then we have the same conclusions as in Proposition \ref{thm1}.
The problem $\Gamma(y_0)$ has a {\it limit value} which is
$V^*(y_0)$. The convergence of $V_t$ to $V^*$  is uniform over
$Z$, and we have $V^*(y_0)$ $=$ $ \sup_{t \geq 1}\; \inf_{m\geq 0}
\; V_{m,t}(y_0)$ $=$ $   \inf_{m\geq 0}\sup_{t\geq 1} \;
V_{m,t}(y_0)$ $=$ $\lim_{m \to \infty, t \to \infty}
V_{m,t}(y_0)$. Moreover the value of $\Gamma(y_0)$  is uniform.
\end{thm}

\noindent {{\bf \underline{Remarks}}:}

$\bullet$ Although $\Delta$ may not satisfy the triangular inequality nor the separation property,
 it   may be seen as a ``distance" adapted to the problem $\Gamma(y_0)$.

$ \bullet$  The assumption: ``$\{(g(y,u), h(y,u)) \in \R^d \times
[0,1] , \;
 u \in U\}$  closed" could be checked for instance if $U$ is
 compact and if $h$ and $g$ are continuous with respect to
 $(y,u)$.

 $\bullet$ $D\uparrow$ is the contingent epi-derivative  (cf \cite{AF}) (which reduces to
 the upper Dini derivative if $\Delta$ is Lipschitz),
  defined by: $D$$\uparrow$$\Delta(z)(\alpha)=
  \liminf_{t \to 0^+, \alpha'\to \alpha} \frac{1}{t}(\Delta(z+t \alpha')-\Delta(z)) $.
  If $\Delta$ is differentiable,   the condition   $D \uparrow \Delta(y_1,y_2)(g(y_1,u),g(y_2,v))\leq 0$
    just reads: $< g(y_1,u), \frac{\partial}{\partial y_1}\Delta (y_1,y_2)> +
    <g(y_2,v),\frac{\partial}{\partial y_2} \Delta (y_1,y_2)>\leq 0$.

$\bullet$ Proposition \ref{thm1} will be a corollary of Theorem
\ref{thm2}. It corresponds to the case where:
$\Delta(y_1,y_2)={\|y_1-y_2\|}^2$, $G(y_0)$ is bounded, and
$h(y,u)=h(y)$ does not depend on $u$ (one can just  take
$\hat{\alpha}(t)=\sup \{|h(x)-h(y)|, {\|x-y\|}^2\leq t\}$).

$\bullet$ $H2a)$  is  a precompacity condition. It is  satisfied as soon as
 $G(y_0)$ is bounded. It is also satisfied if
 $\Delta$ satisfies the triangular inequality and the usual
  precompacity condition: for each $\varepsilon >0,$
  {there exists a finite subset } $C$   {of}  $Z$    {s.t.} : $\forall z \in Z, \exists c \in C,  \Delta(z,c)\leq \varepsilon$. (see lemma \ref{lem4})

$\bullet$  Notice that $H2$ is   satisfied  with  $\Delta=0$ if we
are in the trivial case where   $\inf_u h(y,u)$ is constant.

$\bullet$  Theorem \ref{thm2} can be applied to example 4, with
$\Delta(y_1,y_2)=\|y_1-y_2\|_1$ ($L^1$-norm). In this example, we
have for each $y_1$, $y_2$ and $u$:
$\Delta(y_1+tg(y_1,u),y_2+tg(y_2,u))\leq \Delta(y_1,y_2)$ as soon
as $t\geq 0$ is small enough.

 \subsection{Proof of   Theorem \ref{thm2}}

 We assume in this section that the hypotheses of
 Theorem \ref{thm2} are satisfied,
 and we may assume without loss of generality  that $\hat{\alpha}$ is non decreasing and upper semicontinuous
   (otherwise we replace $\hat{\alpha}(t)$ by $\inf_{\varepsilon >0} \sup_{t' \in [0,t+\varepsilon]}\alpha(t')$).

 \subsubsection{A non expansion property}

 We start with a proposition expressing the fact
 that the problem is non expansive with respect to $\Delta$,
  the idea being that given two initial conditions $y_1$ and $y_2$ and a control
    to be played  at $y_1$,
there exists another control to be played at $y_2$ such that $t
\mapsto  \Delta(y(t,u,y_1),y(t,v,y_2))$ will not increase.

 \begin{pro} \label{pro1} We suppose  the hypothesis of Theorem
 \ref{thm2}. Then
 \be \label{viabcond} \left\{ \begin{array}{l}
  \forall (y_1,y_2)  \in \bar{Z}^2 , \;  \forall T \geq 0, \; \forall  \varepsilon>0, \forall   u \in {\cal U},
  \; \exists v \in {\cal U},\\
\forall t \in [0,T], \;
 \Delta(y(t,u,y_1) , y(t,v,y_2)) \leq
\Delta(y_1,y_2)+ \varepsilon, \\ \mbox{ and for almost every $t
\in
[0,T]$, } \\
h(y(t,v,y_2),v(t))- h(y(t,u,y_1)
  ,u(t)) \leq \hat{\alpha}(\Delta(y(t,u,y_1),y(t,v,y_2) )).
\end{array} \right.\ee
 \end{pro}
\dok First fix $y_1,y_2$ $ \eps >0 $, $T >0$ and $u$. Let us
consider the following set-valued map $\Phi : \R _+ \times \bar{Z}
\times \bar{Z} \times \R \to \R ^d \times \R ^d \times \R $
$$ \Phi (t, x,y,l) := co \, cl  \{ (g(x,u(t) ) , g(y,v) ,0) ) \, | \; v
\in U, \, h(y,v) - h(x,u(t))  \leq \hat{\alpha } (\Delta
(x,y))\},$$ where $co$ stands for the  convex hull and $cl$
for the closure. Notice that $ \Phi (t, x,y,l)$ does not depend on $l$.  Using (\ref{HYPO}), {\it H1)} and {\it H2)b)},  one can check
that $\Phi$ is a set valued map which is upper semicontinuous  in
$(x,y,l)$, measurable in $t$ and with compact convex nonempty
values \cite{AF, Deimling}. We also denote $ \tilde{\Phi} $ the set
valued map defined as $\Phi$ but removing the convex hull.

>From the measurable Viability Theorem \cite{3pol} (cf also
\cite{Vrabie} section 6.5), condition (H2) b) implies that the
epigraph of $\Delta$ (restricted to $\bar{Z}^2 \times
\R$) is viable for the differential inclusion \be \label{diffincl}
(x'(t),y'(t),l'(t)) \in \Phi (t,x(t),y(t),l(t)) \mbox{ for a. e. $
t \geq 0$ }\ee So starting from $(y_1,y_2, \Delta (y_1, y_2) )$,
there exists a solution $(x(\cdot), y(\cdot) , l(\cdot) ) $  to
(\ref{diffincl}) which stays for any $t \geq 0$ in the epigraph of
$\Delta$ namely \be \label{vvv} \Delta (x(t),y(t)) \leq l(t) =
\Delta (y_1,y_2) , \; \forall t \geq 0, \ee by noticing that
$l(\cdot) $ is a constant.

>From the suppositions made on the dynamics $g$, the trajectory
$(x(\cdot), y(\cdot)  ) $ remains in a compact set (included in
some  large enough ball $B(0,M)$)  on the time interval $[0,T]$.
Because $\Delta$ is uniformly continuous on $B(0,M)\times B(0,M)$, there exists $\eta \in
(0,1)  $ with
$$ \forall (x,x',y,y') \in B(0, M+1) ^4,\;  \|x-x'\| + \|y-y'\| <\eta \Longrightarrow |\Delta(x,y) - \Delta
(x',y') | < \eps.$$

Thanks to the Wazewski Relaxation Theorem  (cf for instance Th.
10.4.4 in  \cite{AF})  applied to $\Phi $, the trajectory
$(x(\cdot), y(\cdot), l(\cdot) )$ could be approximated on every
compact interval  by a trajectory to the differential inclusion
defined by $\tilde{\Phi}$. So there exists $(y_1(\cdot) , y_2(\cdot)
, l(\cdot) )$ satisfying
$$ (y_1'(t),y_2'(t),l'(t)) \in \tilde{\Phi} (t,y_1(t),y_2(t),l(t))
\mbox{ for a. e. $ t \geq 0$ }$$ such that $$ \|x(t)-y_1(t)\| +
\|y(t)-y_2(t)\| <\eta ,\; \forall t \in [0,T].$$ From the choice
of $\eta$ and the very definition of $\tilde{\Phi}$ we also  have
    for any
$t \in [0,T]$
\[\left\{
\begin{array} {l} \Delta( y_1(t),y_2(t) ) \leq \Delta (x(t),y(t)) +
\eps \leq \Delta (y_1,y_2) + \eps \\   h(y_2(t),v(t)) - h(y_1(t),u(t))
\leq \tilde{\alpha} (\Delta(y_1(t),y_2(t)) ) \end{array}
\right.  \] This completes our proof if, from one hand we observe
that $y_1 (\cdot) = y(\cdot, u, y_1)$ and from the other hand, we
use Filippov's measurable selection Theorem (e.g. Theorem 8.2.10
in  \cite{AF})  to $\tilde{\Phi}$ for finding a measurable control
$v \in {\cal U} $ such that $y_2(\cdot) = y(\cdot, v ,y_2)$. \QED

\subsubsection{The limit value exists}

Since $\hat{\alpha}$ is u.s.c. and non decreasing, we obtain the
following consequence of  Proposition \ref{pro1}.

\begin{cor} \label{cor1}For every $y_1$ and $y_2$ in $G(y_0)$,   for all $T>0$,
 $$|V_T(y_1)-V_T(y_2)|\leq \hat{\alpha}(\Delta(y_1,y_2)).$$\end{cor}
  Define now, for each $m\geq 0$, $G^m(y_0)$ as the set of states which can be reached from $x_0$ before time $m$:
  $$G^m(y_0)=\{y(t,u,y_0) , t\leq m, u \in {\cal U}\}, \; \; {\rm so  \; that\;\;} G(y_0)= \cup_{m \geq 0}G^m(y_0).$$
     An immediate consequence of the precompacity hypothesis H2a) is the following
\begin{lem} \label{lem2} For every $\varepsilon>0$,  there exists $m_0$ in $\R_+$ such that:
 {$$\forall z \in G(y_0),  \exists z' \in G^{m_0}(y_0)\;  \mbox{ such that }   \; \Delta(z,z')\leq \varepsilon.$$}
\end{lem}
\dok  Otherwise for each positive integer $m$ one can find $z_m$
in $G(y_0)$ such that $\Delta(z_m,z)>\varepsilon$ for all $z$ in
$G^m(y_0)$. Use {\it H2a)}   to find  $n$ such that $\liminf_m
\Delta(z_n,z_m)\leq \varepsilon$. Since $z_n\in G(y_0)$, there
must exist $k$ such that $z_n\in G^k(y_0)$. But for each $m\geq k$
we have $z_n \in G^m(y_0)$, hence  $\Delta(z_m, z_n)>\varepsilon$.
We obtain  a contradiction. \QED

\vspace{0,3cm}

We can already conclude for the limit value.

 \begin{pro} \label{pro2}  $V_t(y_0) \xrightarrow[t \to \infty]{} V^*(y_0).$
 \end{pro}

 \dok  Because of lemma \ref{lem1},   it is sufficient to prove that for
 every $\varepsilon>0$, there exists $m_0$ such that:
$$\sup_{t>0} \inf_{m\leq m_0} V_{m,t}( y_0) \leq    \sup_{t>0} \inf_{m\geq 0}   V_{m,t}(y_0) + 2 \eps $$
Fix $\varepsilon$, and consider $\eta>0$ such that $\hat{\alpha}(t)\leq \varepsilon$ as soon as
 $t\leq \eta$. Use lemma \ref{lem2} to find $m_0$ such that $\forall z \in G(y_0),  \exists z' \in G^{m_0}(y_0)\;
  s.t.  \; \Delta(z,z')\leq \eta.$

Consider any $t>0$. We   have  $ \inf_{m\geq 0}  V_{m,t}( y_0) =
\inf \{ V_t(z), z \in G(y_0)\}$, and $ \inf_{m\leq m_0}   V_{m,t}(
y_0) = \inf \{ V_t(z), z \in G^{m_0}(y_0)\}$. Let $z$ in $G(y_0)$
be such that $V_t(z)\leq  \inf_m V_{m,t}( y_0)+ \varepsilon$, and
consider $z' \in G^{m_0}(y_0)\;  s.t.  \; \Delta(z,z')\leq \eta$.
By corollary \ref{cor1},  $|V_t(z)-V_t(z')|\leq
\hat{\alpha}(\Delta(z,z')) \leq \eps ,$ so we obtain that
$$\inf_{m\leq m_0} V_{m,t}( y_0) \leq V_t(z') \leq V_t(z)+ \eps
\leq \inf_{m} V_{m,t}( y_0) + 2 \eps.$$ Passing to the
supremum on $t$, this completes the proof. \QED

\begin{rem}\label{rq1} Observe that for obtaining the existence of the value, we have used a compactness argument
(assumption H2)a)) and condition (\ref{viabcond}). We did not use
explicitly assumption H2)b)  which is only used for obtaining
(\ref{viabcond}).
\end{rem}

The rest  of the proof is more involved, and is  inspired by the
proof of Theorem 3.6 in \cite{Renault}.

\subsubsection{Auxiliary value functions}

The uniform value requires the same control to be good for all  time horizons,
and we are led to  introduce new auxiliary value functions.
 Given $m\geq 0$ and $n\geq 1$,  for any initial state $z$ in $Z=G(y_0)$ and control
 $u$ in ${\cal U}$, we define
$$\nu_{m,n}(z,u)=\sup_{t \in [1,n]} \gamma_{m,t}(z,u),\; {\rm and }\;
 W_{m,n}(z)=  \inf_{u \in {\cal U}} \nu_{m,n}(z,u).$$

\noindent $W_{m,n}$ is the value function of the problem where the controller can
 use the time interval $[0,m]$ to reach a good state, and
 then his cost is only the supremum for $t$ in $[1,n]$, of the average cost between time $m$ and $m+t$.
 Of course, we have $W_{m,n}\geq V_{m,n}$. We write $\nu_{n}$ for $\nu_{0,n}$, and $W_n$ for $W_{0,n}.$

We easily obtain from proposition \ref{pro1}, as in corollary \ref{cor1}, the following result.

\begin{lem} \label{lem2,5}For every $z$ and $z'$ in $Z$,   for all $m\geq 0$ and $n\geq 1$,
 $$|V_{m,n}(z)-V_{m,n}(z')|\leq \hat{\alpha}(\Delta(z,z')).$$
 $$|W_{m,n}(z)-W_{m,n}(z')|\leq \hat{\alpha}(\Delta(z,z')).$$\end{lem}

 The following lemma   shows that the quantities $W_{m,n}$ are not that high.
 \begin{lem}  \label{lem3}$\forall k \geq 1,
 \forall n \geq 1, \forall m \geq 0,\forall z \in Z,$  $$V_{m,n}(z) \geq
 \inf_{l \geq m} W_{l,k}(z) - \frac{k}{n}.$$\end{lem}
\dok Fix $k$, $n$, $m$ and $z$, and put $A= \inf_{l \geq m}
W_{l,k}(z)$. Consider any control $u$ in ${\cal U}$.
  For any $i\geq m$, we have
   $$\sup_{t \in [1, k]}  \gamma_{i,t}(z,u)= \nu_{i,k}(z,u) \geq W_{i,k}(z)\geq A.$$ So
    we know that for any $i\geq m$, there exists $t(i)\in [1,k]$ such that  $\gamma_{i,t(i)}(z,u) \geq A$.

 Define now by induction $i_1=m$, $i_2=i_1+t(i_1)$,..., $i_q=i_{q-1}+t(i_{q-1})$,
  where $q$ is such that $i_q\leq n +m < i_q +t(i_q)$.
  We have $n \gamma_{m,n}(z,u) \geq \sum_{p=1}^{q-1}t(i_p) A  \geq nA- k$,
   so $\gamma_{m,n}(z,u) \geq A - \frac{k}{n}.$ Taking the infimum
   over all controls, the proof is complete. \QED

 \vspace{0,3cm}

 We know from Proposition \ref{pro2} that the limit value is given by $V^*$. We now give other  formulas for this limit.
  \begin{pro}  \label{pro2,1}For every state $z$ in $Z$,
$$ \inf_{m \geq 0}   \sup _{n\geq 1}\;\; W_{m,n}(z)=
 \inf_{m \geq 0}   \sup_{n\geq 1}\;\; V_{m,n}(z)=V^*(z) = \sup_{n \geq 1} \inf _{m \geq 0} V_{m,n}(z)
=\sup_{n\geq 1} \inf_{m \geq 0}  W_{m,n}(z).$$ \end{pro}

  \noindent {\bf Proof of proposition  \ref{pro2,1}:} Fix an initial state $z$ in $Z$.
  We already  have   $V^*(z)= \sup_{t> 0}\; \inf_{m\geq 0} \; V_{m,t}(z)\geq \sup_{t\geq 1}\;
   \inf_{m\geq 0} \; V_{m,t}(z)$.
 One can easily check that $\inf_{m\geq 0} \; V_{m,t}(z)\leq \inf_{m\geq 0} \; V_{m,2t}(z)$ for each positive
 $t$. So
 $$V^*(z) \geq \sup_{t\geq 1}\;
   \inf_{m\geq 0} \; V_{m,t}(z) \geq \sup_{t\geq (1/2)}\;
   \inf_{m\geq 0} \; V_{m,t}(z) \geq \ldots \sup_{t> 0}\; \inf_{m\geq 0} \;
   V_{m,t}(z) =V^*(z).$$
Consequently $V^*(z) =\sup_{t\geq 1}\;
   \inf_{m\geq 0} \; V_{m,t}(z)$. Moreover because $V_{m,t} \leq
   W_{m,t}$ we have also $V^*(z) \leq \sup_{t\geq 1}\;
   \inf_{m\geq 0} \; W_{m,t}(z).$

  We now claim   that  $V^*(z) = \sup_{t\geq 1}\;
   \inf_{m\geq 0} \; W_{m,t}(z)$. It remains to show
 $V^*(z) \geq \sup_{t\geq 1}\;
   \inf_{m\geq 0} \; W_{m,t}(z)$. From Lemma \ref{lem3}, we know
   that for all
 $k \geq 1$, $ n \geq1$ and  $ m \geq 0$, we have
  $V_{m,nk}(z)\geq  \inf_{l \geq 0} W_{l,k}(z) - \frac{1}{n}$, so
  $\inf_m V_{m,nk}(z) \geq \inf_{l\geq 0} W_{l,k}(z) -\frac{1}{n}$.
By taking the supremum on $n$ , we obtain
 $$ V^ * (z) =\sup_{n\geq 1} \inf_{m\geq 0} V_{m,n}(z) \geq \sup_{n \geq 1} \inf_{m \geq 0}  V_{m,nk}(z)
  \geq \inf_{l \geq 0} W_{l,k}(z).$$ Since $k$ is arbitrary, we have proved our claim.

  Since the inequalities
  $$ \inf_{m \geq 0}   \sup_{n\geq 1}\;\; W_{m,n}(z)\geq \inf_{m \geq 0}   \sup_{n\geq 1}\;\; V_{m,n}(z)
  \geq  \sup_{n\geq 1}\;\;  \inf_{m \geq 0}  V_{m,n}(z)= V^*(z)$$ are clear, to conclude the proof of the proposition
   it is enough to  show that $ \inf_{m \geq 0}   \sup_{n\geq 1}\;\; W_{m,n}(z) \leq V^*(z)$.

  Fix $\varepsilon >0$. We have already proved that  $V^*(z)=\sup_{n\geq 1} \inf_{m \geq 0}  W_{m,n}(z)$,
   so for each $n \geq 1$ there exists $m\geq 0$ such that $W_{m,n}(z)\leq V^*(z)+\varepsilon.$
   Hence for each $n$, there exists $z'_n$ in $G(z)$ such that $ W_{0,n}(z'_n)
   \leq V^*(z) + \varepsilon $. We know from Lemma \ref{lem2} that  there exists $m_0\geq 0$ such that:
 {$\forall z' \in G(z),  \exists z'' \in G^{m_0}(z)\;  s.t.  \; \Delta(z',z'')\leq \varepsilon.$}
 Consequently, for each $n\geq 1$, there exists $z''_n$ in $G^{m_0}(z)$ such that
 $ \Delta(z'_n,z''_n)\leq \varepsilon$, and by lemma \ref{lem2,5} this implies that
  $$W_n(z''_n) \leq W_n(z'_n) + \hat{\alpha}(\varepsilon)\leq V^*(z)+  \varepsilon +
  \hat{\alpha}(\varepsilon).$$

 Up to now, we have proved that for every $\varepsilon ' >0$, there exists  $m_0$ such that:
 $$\forall n \geq 1, \exists m \leq m_0 \; {\rm s.t.} \; W_{m,n}(z)\leq V^*(z) + \varepsilon '.$$
 Since all costs lie in $[0,1]$, it is easy to check that $|W_{m,n}(z)-W_{m',n}(z)| \leq |m-m'|$ for each
  $n$,  $m$, $m'$. Hence there exists a {\it finite} subset $F$ of $[0,m_0]$ such that:
  $\forall n \geq 1, \exists m \in F  \; {\rm s.t.} \; W_{m,n}(z)\leq V^*(z)+ 2 \varepsilon '.$
   Considering $\hat{m}$ in $F$ such that the set
   $\{n \; {\rm positive \;  integer,}\; W_{\hat{m},n}(z)\leq V^*(z)+ 2 \varepsilon ' \}$
    is infinite, and noticing that
 $W_{m,n}$ is non decreasing in $n$, we obtain the existence of a unique
  $\hat{m}\geq 0$ such that  $\forall n \geq 1,     \; W_{\hat{m},n}(z)\leq V^*(z)+2 \varepsilon'.$
   Hence $ \eps '$ being arbitrary,  $ \inf_{m \geq 0}   \sup_{n\geq 1}\;\; W_{m,n}(z) \leq V^*(z)$, concluding the proof
   of Proposition \ref{pro2,1}. \QED

  \vspace{0,3cm}

   We now look for uniform convergence properties. By the precompacity condition $H2a)$, it is easy to obtain that:

\begin{lem} \label{lem4}  For each $\varepsilon >0,$  {there exists a finite subset } $C$   {of}  $Z$    {s.t.} : $\forall z \in Z, \exists c \in C,  \Delta(z,c)\leq \varepsilon. $ \end{lem}

We know that $(V_n)_n$ simply converges to $V^*$ on $Z$.  Since $|V_{n}(z)-V_{n}(z')|\leq \hat{\alpha}(\Delta(z,z'))$ for all $n$, $z$ and $z'$, we obtain by lemma \ref{lem4}:

\begin{cor} \label{cor2 }  The convergence of $(V_n)_n$ to $V^*$ is uniform on $Z$. \end{cor}

\vspace{0,5cm}

We can proceed similarly to obtain other uniform properties. We
have $$\displaystyle{V^*(z)=\sup_{n \geq 1} \inf_{m \geq 0}
W_{m,n}(z) = \lim _{n \to + \infty } \inf_{m \geq 0} W_{m,n}(z) }
$$  since $\inf_{m \geq 0}  W_{m,n}(z)$ is not decreasing in
$n$. Using lemmas \ref{lem2,5} and \ref{lem4}, we obtain that the
convergence is uniform, hence we get:
   $$\forall \varepsilon>0, \exists n_0, \forall z \in Z, \;\;
   V ^ * (z) - \eps  \leq \inf_{m\geq0} W_{m,n_0} (z) \leq V^*(z).$$
By Lemma \ref{lem3}, we obtain:$$\forall \varepsilon>0, \exists
n_0, \forall z \in Z, \forall m\geq 0, \forall n \geq 1,
V_{m,n}(z) \geq V^*(z) - \varepsilon - \frac{n_0}{n}.$$
Considering $n$ large  gives:
\begin{equation}\label{eq3} \forall \varepsilon>0,
\exists K, \forall z \in Z, \forall n \geq K,  \;\; \inf_{m\geq 0}
V_{m,n}(z) \geq V^*(z)- \varepsilon\end{equation}

 Write now, for each state $z$ and $m\geq 0$: $h_m(z)= \inf_{m'\leq m}  \sup_{n\geq 1} W_{m',n}(z)$.
 $(h_m)_m$   converges to   $V^*$, and as before, by Lemmas
  \ref{lem2,5} and \ref{lem4}, we obtain that the convergence is uniform. Consequently,
 \begin{equation} \label{eq4}  \forall \varepsilon>0, \exists M\geq 0,
  \forall z \in Z, \exists m\leq M, \;\; \sup_{n \geq 1} W_{m,n}(z) \leq V^*(z)+ \varepsilon.\end{equation}

 \subsubsection{On the existence of a uniform value }

  In order to prove that $\Gamma(y_0)$
   has a uniform value we have
   to show that for every $\varepsilon>0$,
   there exist a control $u$ and a time $n_0$ such that for every $n\geq n_0$,
    $\gamma_{n}(y_0,u)\leq V^*(y_0)+\varepsilon$.
     In this subsection we  adapt the proofs of Lemma 4.1 and Proposition 4.2.
     in \cite{Renault}.
 We start by constructing, for each $n$, a control which: 1)
  gives  low average costs if one stops the play at
  {\it any} large time  before $n$, and 2) after time $n$,
  leaves the player with a good ``target" cost.
   This explains  the importance of the quantities $\nu_{m,n}$. We start with
   the following
\begin{lem} \label{lem6}
$\forall \varepsilon>0, \exists M \geq0,  \exists K \geq 1,
\forall z \in Z, \exists m \leq M, \forall n \geq K, \exists u \in
{\cal U}$ such that: \be \label{qqq} \nu_{m,n } (z,u)\leq V^*(z)+
\varepsilon/2, \; \mathnormal{\rm and }\;  V^*(y(m+n, u,z) )\leq
V^*(z)+ \varepsilon.\ee
\end{lem}

\dok Fix $\varepsilon>0$. Take $M$ given by (\ref{eq4}), so that $
\forall z \in Z, \exists m\leq M, \;\; \sup_{n \geq 1} W_{m,n}(z)
\leq V^*(z)+ \varepsilon.$ Take $K\geq 1$ given by (\ref{eq3})
such that: $\forall z \in Z$, $\forall n \geq K$, $\inf_{m}
V_{m,n}(z) \geq V^*(z) - \varepsilon.$

Fix an initial state $z$ in $Z$. Consider $m$ given by
(\ref{eq4}), and $n \geq K.$ We have to find $u$ in ${\cal U}$
satisfying (\ref{qqq}).

We have $W_{m,n'}(z) \leq V^*(z)+ \varepsilon$ for every $n'\geq
1$, so $W_{m,2n}(z)\leq V^*(z)+\varepsilon$, and we consider a
control $u$  which is $\varepsilon$-optimal for $W_{m,2n}(z)$, in
the sense that $\nu_{m,2n}(z,u)\leq W_{m,2n}(z)+\varepsilon.$ We
have:

\centerline{$\nu_{m,n}(z,u) \leq \nu_{m,2n}(z,u) \leq W_{m,2n}(z)+
\varepsilon \leq V^*(z)+ 2\varepsilon.$}

\noindent Denoting $X=\gamma_{m,n}(z,u)$ and
$Y=\gamma_{m+n,n}(z,u)$.
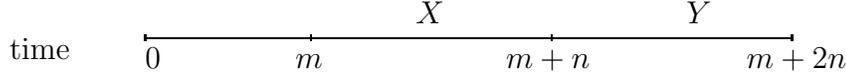
\begin{center} \setlength{\unitlength}{0,4mm}
\begin{picture}(230,10)
\put(0,0){\line(1,0){55}}
\put(55,0){\line(1,0){80}}
\put(135,0){\line(1,0){80}}

\put(0,-1){\line(0,1){2}}

\put(55,-1){\line(0,1){2}}
\put(135,-1){\line(0,1){2}}

\put(215,-1){\line(0,1){2}}

 \put(-45,-7){${\rm time}$}
\put(0,-10){$0$}

\put(50,-10){$m$}
\put(120,-10){$m+n$}

 \put(200,-10){${m+2n}$}

\put(90,5){$X$}
\put(180,5){$Y$}
\end{picture}

      \end{center}

\vspace{0,3cm}

 \noindent Since $\nu_{m,2n}(z,u)\leq V^*(z)+2\varepsilon$, we have
 $X\leq V^*(z)+2\varepsilon$, and $(X+Y)/2= \gamma_{m,2n}(z,u)\leq V^*(z)+2\varepsilon.$
  Since $n\geq K$, we also have $X\geq V_{m,n}(z) \geq V^*(z)- \varepsilon$.
  And $n \geq K$ also gives $V_n(y(m+n, u,z)) \geq V^*(y(m+n, u,z ))-  \varepsilon$,
  so $V^*(y(m+n, u,z))\leq V_n(y(m+n, u,z)) +\varepsilon \leq Y+ \varepsilon.$
  Writing  now  $Y/2=(X+Y)/2-X/2$ we obtain  $Y/2\leq (V^*(z)+ 5 \varepsilon)/2$.
   So $Y\leq V^*(z)+5 \varepsilon$, and  finally $V^*(y(m+n, u,z) )\leq V^*(z)+ 6 \varepsilon.$
   \QED

\vspace{1cm}

 We can now conclude the proof of theorem \ref{thm2}.

  \begin{pro} \label{pro3} For every state $z$ in $Z$ and $\varepsilon>0$, there exists a control $u$
  in ${\cal U}$ and $T_0$ such that for every $T\geq T_0$, $\gamma_{T}(z,u)\leq V^*(z)+\varepsilon$.  \end{pro}

\dok  Fix $\alpha>0$.

For every positive integer  $i$, put $\varepsilon_i=\frac{\alpha}{2^i}$.
Define $M_i=M(\varepsilon_i)$ and $K_i=K(\varepsilon_i)$ given by lemma \ref{lem6} for $\varepsilon_i$.
Define also $n_i= \Max \{K_i, \frac{M_{i+1}}{\alpha}\}\geq 1$.

 We have: $ \forall i \geq 1, \forall z \in Z, \exists \,
m(z,i)\leq M_i, \exists u \in {\cal U}, \; \rm s.t.$
$$\; \nu_{m(z,i),n_i}(z,u) \leq V^*(z)+\frac{\alpha}{2^{i+1}}\; {\rm and }\;  V^*(y(m(z,i)+n_i,u,z))
\leq V^*(z)+\frac{\alpha}{2^i}.$$

We now fix the initial state $z$ in $Z$, and for simplicity write $v^*$ for $V^*(z)$.
  We define a sequence $(z^i,m_i,u^i)_{i \geq 1}$ by induction:

$\bullet$ first put $z^1=z$, $m_1=m(z^1,1)\leq M_1$,
 and pick $u^1$ in  ${\cal U}$ such that
 $\nu_{m_1,n_1}(z^1,u^1)\leq V^*(z^1)+\frac{\alpha}{2^2}$, and $V^*(y(m_1+n_1, u^1,z^1) )
 \leq V^*(z^1)+ \frac{\alpha}{2}.$

$\bullet$ for $i\geq 2$, put $z^i= y(m_{i-1}+n_{i-1}, u^{i-1} ,  z^{i-1})$,
 $m_i=m(z^i, i)\leq M_i$, and pick $u^i$ in
  ${\cal U}$  such that $\nu_{m_i,n_i}(z^i,u^i)\leq V^*(z^i)+\frac{\alpha}{2^{i+1}}$ and
   $V^*(y( m_i +n_i, u^{i}, z^{i}) )\leq V^*(z^i)+\frac{\alpha}{2^i}.$\\

Consider finally $u$ in ${\cal U}$ defined by concatenation:
 first $u^1$ is followed for time $t$ in $[0, m_1+n_1)$, then $u^2$ is followed for
 $t$ in $[m_1+n_1, m_2+n_2)$, etc... Since
 $z^i= y(m_{i-1}+n_{i-1}, u^{i-1} ,  z^{i-1})$ for each $i$,
 we have $y(\sum_{j=1}^{i-1} m_{j}+n_{j},u,z)= z^i$ for each $i$.  For each $i$
  we have $n_i \geq M_{i+1}/\alpha \geq m_{i+1}/ \alpha$, so an
   interval with length $n_i$   is much longer than an interval with length $m_{i+1}$.

 \begin{center}
\setlength{\unitlength}{0,4mm}
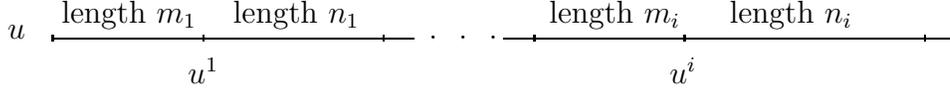
\begin{picture}(300,15)
\put(0,0){\line(1,0){120}}
\put(150,0){\line(1,0){150}}
\put(125,0){.}
\put(135,0){.}
\put(145,0){.}

\put(0,-1){\line(0,1){2}}
\put(50,-1){\line(0,1){2}}
\put(110,-1){\line(0,1){2}}
\put(160,-1){\line(0,1){2}}
\put(210,-1){\line(0,1){2}}
\put(290,-1){\line(0,1){2}}

 \put(-15,0){$u$}

\put(45,-15){$u^1$}
\put(205,-15){$u^i$}

\put(3,5){length $m_1$}
\put(60,5){length $n_1$}
\put(165,5){length $m_i$}
\put(225,5){length $n_i$}
\end{picture}

      \end{center}

\vspace{0,5cm}

For each $i\geq 1$, we have $V^*(z^{i})\leq
 V^*(z^{i-1})+\frac{\alpha}{2^{i-1}}$.
  So $V^*(z^{i})\leq +\frac{\alpha}{2^{i-1}} +\frac{\alpha}{2^{i-2}}...  +\frac{\alpha}{2} + V^*(z^1)\leq v^* +\alpha
  -\frac{\alpha}{2^i}.$ So $\nu_{m_i,n_i}(z^i,u^i)\leq v^* + \alpha$. \\

Let now $T$ be large.

- First assume that $T=m_1+n_1+...+m_{i-1}+n_{i-1}+r$, for some positive $i$ and $r$ in $[0,m_i]$. We have:
 \begin{eqnarray*}
 \gamma_T(z,u) & = &\frac{1}{T}   \int_{0}^{T} h(y(s,u,z),u(s))ds \\
 \; & \leq &  \frac{1}{T}  \left(\sum_{j=1}^{i-1} n_j\right)
 (v^*+\alpha)+ \frac{m_1}{T} + \frac{1}{T}  \left(\sum_{j=2}^{i} m_j\right)
 \end{eqnarray*}
But $m_j\leq \alpha n_{j-1}$ for each $j$, so
 $$ \gamma_T(z,u)\leq v^*+ 2\alpha+ \frac{m_1}{T}.$$

-  Assume now that $T=m_1+n_1+...+m_{i-1}+n_{i-1}+m_i+r$, for some positive $i$ and $r$ in $[0,n_i]$.
 The previous computation shows that:
 $$ \int_{0}^{T-r} h(y(s,u,z),u(s))ds\leq m_1 +(T-r) (v^*+ 2\alpha).$$ Since
 $\nu_{m_i,n_i}(z^i,u^i)\leq v^* + \alpha$, we obtain:
  \begin{eqnarray*}
  T \gamma_T(z,u) & = & \int_0^{T-r} h(y(s,u,z) ,u(s))ds  +  \int_{T-r}^T h(y(s,u,z) ,u(s))ds, \\
  & \leq &  m_1 +(T-r) (v^*+ 2\alpha) + r (v^* + \alpha),\\
 & \leq &  m_1 +T  (v^*+ 2\alpha).
 \end{eqnarray*}
  Consequently, here also we have:  $$ \gamma_T(z,u)\leq v^*+ 2\alpha+ \frac{m_1}{T}.$$

  This concludes the proofs of Proposition \ref{pro3} and
consequently, of Theorem \ref{thm2}. \QED

 \vspace{1cm}

\noindent \large \bf Acknowledgements. \rm \small

The first author wishes to thank Pierre Cardaliaguet, Catherine Rainer and Vladimir Veliov for stimulating conversations. The second author wishes to thank Patrick Bernard, Pierre Cardaliaguet,   Antoine Girard,   Filippo Santambrogio and Eric S\'er\'e for fruitful discussions.

The work of Jerome Renault was partly supported by the   French
Agence Nationale de la Recherche (ANR), undergrants  ATLAS     and
Croyances, and the ``Chaire de la Fondation du Risque",
Dauphine-ENSAE-Groupama : Les particuliers face aux risques.


\vspace{1cm}




\begin{thebibliography}{9}




\bibitem{A2} Arisawa, M. and P.L. Lions    (1998)
\newblock Ergodic problem for the Hamilton Jacobi Belmann equations II,
\newblock {\em  Ann. Inst. Henri Poincar\'e, Analyse Nonlin\'eaire}, 15 ,1 ,  1--24.

\bibitem{Arisawa} Arisawa, M. and P.L. Lions    (1998)
\newblock On ergodic stochastic control.
\newblock {\em Com. in partial differential equations}, 23, 2187--2217.

\bibitem{artstein} Z. Artstein, and V. Gaitsgory,
{\em The value function of singularly perturbed control
systems\/}, Appl. Math. Optim., 41 (2000), 425-445.

\bibitem{AC} Aubin J. P.,   A. Cellina   (1984)
\newblock Differential Inclusion
\newblock {\em Springer}.

\bibitem{AV} Aubin J. P.,     (1992) \newblock Viability Theory \newblock {\em Birkhauser}.

\bibitem{AF} Aubin J. P.,  Frankowska H.   (1990)
\newblock Set-Valued Analysis
\newblock {\em Birkha\"user}.


\bibitem{Bettiol} Bettiol, P. (2005)
\newblock On ergodic problem for Hamilton-Jacobi-Isaacs equations
\newblock {\em ESAIM: Cocv}, 11, 522--541.

\bibitem{Carda} Cardaliaguet P.
Ergodicity of Hamilton-Jacobi equations with a non coercive non
convex Hamiltonian in $\R^2/ Z^2$  preprint [hal-00348219 -
version 1] (18/12/2008)

\bibitem{Vrabie} Carja, O., Necula, M., Vrabie, I. (2007)  Viability,
Invariance and Applications, North-Holland.

\bibitem{Deimling} K. Deimling (1992) Multivalued Differential
Equations, De gruyter Seris in Nonlinear Analysis and
Apllications.

\bibitem{3pol} Frankowska, H., Plaskacz, S. and Rzezuchowski T. (1995):
Measurable Viability Theorems and Hamilton-Jacobi-Bellman
Equation, J. Diff. Eqs., 116, 265-305.





\bibitem{Lions}
Lions P.-L. , Papanicolaou G. , Varadhan S.R.S., Homogenization of
Hamilton- Jacobi Equations, unpublished work.



\bibitem{Renault}  Renault, J. (2007) Uniform value in Dynamic Programming.  Cahier du Ceremade 2007-1. arXiv : 0803.2758.

\bibitem{Ti} Tichonov, A. N. (1952)
Systems of differential equations containing small parameter near
derivatives, Math. Sbornik. 31  575ñ586.


\bibitem{Veliov} Veliov, V.  Critical values in long time
optimal control. Unpublished work. Seminar of Applied Analysis in
universit\'e de Brest (2003).
\end{thebibliography}
\end{document}